\documentclass{amsart}
\begin{document}

\vfuzz2pt 
\hfuzz2pt 
\newtheorem{thm}{Theorem}[section]
\newtheorem{corollary}[thm]{Corollary}
\newtheorem{lemma}[thm]{Lemma}
\newtheorem{proposition}[thm]{Proposition}
\newtheorem{defn}[thm]{Definition}
\newtheorem{remark}[thm]{Remark}
\newtheorem{example}[thm]{Example}
\newtheorem{fact}[thm]{Fact}
\
\newcommand{\norm}[1]{\left\Vert#1\right\Vert}
\newcommand{\abs}[1]{\left\vert#1\right\vert}
\newcommand{\set}[1]{\left\{#1\right\}}
\newcommand{\Real}{\mathbb R}
\newcommand{\eps}{\varepsilon}
\newcommand{\To}{\longrightarrow}
\newcommand{\BX}{\mathbf{B}(X)}
\newcommand{\A}{\mathcal{A}}
\newcommand{\onabla}{\overline{\nabla}}
\newcommand{\hnabla}{\hat{\nabla}}


\def\proof{\medskip Proof.\ }
\font\lasek=lasy10 \chardef\kwadrat="32 
\def\kwadracik{{\lasek\kwadrat}}
\def\koniec{\hfill\lower 2pt\hbox{\kwadracik}\medskip}

\newcommand*{\C}{\mathbf{C}}
\newcommand*{\R}{\mathbf{R}}
\newcommand*{\Z}{\mathbf {Z}}

\def\sb{f:M\longrightarrow \C ^n}
\def\det{\hbox{\rm det}\, }
\def\detc{\hbox{\rm det }_{\C}}
\def\i{\hbox{\rm i}}
\def\tr{\hbox{\rm tr}\, }
\def\rk{\hbox{\rm rk}\,}
\def\vol{\hbox{\rm vol}\,}
\def\Im {\hbox{\rm Im}\, }
\def\Re{\hbox{\rm Re}\, }
\def\interior{\hbox{\rm int}\, }
\def\e{\hbox{\rm e}}
\def\pai{\partial _i}
\def\paj{\partial _j}
\def\pu{\partial _u}
\def\pv{\partial _v}
\def\pui{\partial _{u_i}}
\def\puj{\partial _{u_j}}
\def\puk{\partial {u_k}}
\def\div{\hbox{\rm div}\,}
\def\lg{\hbox{\rm ln}\,\rho}
\def\Ric{\hbox{\rm Ric}\,}
\def\r#1{(\ref{#1})}
\def\ker{\hbox{\rm ker}\,}
\def\im{\hbox{\rm im}\, }
\def\I{\hbox{\rm I}\,}
\def\id{\hbox{\rm id}\,}
\def\exp{\hbox{{\rm exp}^{\tilde\nabla}}\.}
\def\cka{{\mathcal C}^{k,a}}
\def\ckplusja{{\mathcal C}^{k+1,a}}
\def\cja{{\mathcal C}^{1,a}}
\def\cda{{\mathcal C}^{2,a}}
\def\cta{{\mathcal C}^{3,a}}
\def\c0a{{\mathcal C}^{0,a}}
\def\f0{{\mathcal F}^{0}}
\def\fnj{{\mathcal F}^{n-1}}
\def\fn{{\mathcal F}^{n}}
\def\fnd{{\mathcal F}^{n-2}}
\def\Hn{{\mathcal H}^n}
\def\Hnj{{\mathcal H}^{n-1}}
\def\emb{\mathcal C^{\infty}_{emb}(M,N)}
\def\M{\mathcal M}
\def\Ef{\mathcal E _f}
\def\Eg{\mathcal E _g}
\def\Nf{\mathcal N _f}
\def\Ng{\mathcal N _g}
\def\Tf{\mathcal T _f}
\def\Tg{\mathcal T _g}
\def\diff{{\mathcal Diff}^{\infty}(M)}
\def\embM{\mathcal C^{\infty}_{emb}(M,M)}
\def\U1f{{\mathcal U}^1 _f}
\def\Uf{{\mathcal U} _f}
\def\Ug{{\mathcal U} _g}
\def\hnu{\hat\nu}
\def\gnu{\nu_g}
\def\C{{\mathcal ()}}
\def\A{{\mathcal (*)}}
\def\T{{\rm T}}
\title{Affine spheres with prescribed Blaschke metric}
\author{Barbara Opozda}




\address{Faculty of Mathematics and Computer Science UJ, ul. {\L}ojasiewicza 6, 30-348, Cracow, Poland}

 \email{Barbara.Opozda@im.uj.edu.pl}
 \maketitle
 \begin{abstract}
It is proved that the equality $\Delta \ln|\kappa-\lambda|=6\kappa$,
where $\kappa$ is the Gaussian curvature of a metric tensor $g$ on a
2-dimensional manifold is a sufficient and necessary condition for
local realizability of the metric as the Blaschke metric of some
affine sphere. Consequently, the set of all improper local affine
spheres with nowhere-vanishing Pick invariant can be parametrized by
harmonic functions.
 \end{abstract}

\section{Introduction} Affine spheres are a still a mysterious
class, even in the 2-dimensional case. The following fact has been
known since Blaschke's times, \cite{B}. For an affine sphere in
$\R^3$ whose Pick invariant is nowhere vanishing the following
equality is satisfied
\begin{equation}\label{wstep1}
\Delta \ln|J|=6\kappa,
\end{equation}
where $\kappa$ is the Gaussian curvature of the Blaschke metric and
$J$ is the Pick invariant. The Laplacian is taken relative to the
Blaschke metric. The affine theorema egregium says that $J=\kappa
-\rho$, where $\rho$ is the affine scalar curvature. For an affine
sphere with the shape operator $S=\lambda\, \id$ the scalar
curvature is equal to $\lambda$. Therefore, the equality
(\ref{wstep1}) can be written as $
\Delta\ln|\kappa-\lambda|=6\kappa. $ It has not been noticed,
however,  that this equality is also a sufficient condition for
local realizability of a metric as the Blaschke metric on an affine
sphere in $\R^3$. Realizability of prescribed objects on
submanifolds belongs to fundamental problems in all types of
geometry of submanifolds. The main aim of this paper is to prove the
following result

\begin{thm}\label{main} Let $g$ be a metric tensor field on a $2$-dimensional
manifold $M$. It can be locally realized as the Blaschke metric on
an affine sphere with nowhere-vanishing Pick invariant if and only
if
\begin{equation}\label{warunek}
\Delta \ln |\kappa-\lambda|=6\kappa
\end{equation}
for some constant $\lambda$ such that $\kappa-\lambda\ne 0$
everywhere.
\end{thm}

Note that the 2-dimensional affine sphere with vanishing Pick
inviariant are classified. Namely, by the affine theorema egregium
we know that if an affine sphere has vanishing Pick invariant then
the Gaussian curvature of the Blaschke metric is constant. A
complete  classification of such affine spheres in $\R^3$ is given,
for instance, in Section 5 of  Chapter III in \cite{NS}.

 As corollaries of Theorem \ref{main} we shall prove

 \begin{corollary}\label{corollary_1}
An improper locally strongly convex  affine sphere in $\R^3$  with
nowhere-vanishing Pick invariant is analytic.
 \end{corollary}
Note that  not all improper affine spheres  are analytic. For
instance, the surface given by the equation $z=xy +\Phi(x)$, where
$\Phi$ a smooth but non-analytic function is an affine sphere. This
sphere has vanishing Pick invariant and  its Blaschke metric is
indefinite, see \cite{NS}.

In  the following corollary $g_0$  stands for the canonical metric
tensor field in $\R^2$ and the harmonicity is relative to $g_0$.
\begin{corollary}\label{corollary_2}
Let $\overline g=e^Lg_0$ be a metric tensor of constant Gaussian
curvature $\pm 2$ defined in a neighborhood $U$ of $p\in\R^2$. Let
$h$ be a harmonic function defined on  $U$. Then
$g=e^{\frac{h-L}{2}}g_0$ can be locally realized as the Blaschke
metric in a neighborhood of $p$   on some improper affine sphere.
\end{corollary}

Since metric structures of the same constant Gaussian curvature are
locally  isometric, the above theorem says, roughly speaking, that
local  improper affine spheres with nowhere-vanishing Pick invariant
can be parametrized by harmonic functions.

\bigskip

\section{Affine spheres}

 Let $f:M\to \R^{n+1}$ be a hypersurface in $\R^{n+1}$. At the beginning, for
 simplicity, we assume that $M$ is connected and oriented.
 On $\R^{n+1}$ we have the volume form given by the standard
 determinant. The standard flat connection on $\R^{n+1}$ will be denoted by $\tilde \nabla$.
 Let $\xi$ be a transversal vector field for $f$ (consistent with the orientation on $M$).
The induced volume form on $M$ is given by
\begin{equation}
\nu _\xi(X_1,...,X_n)=\det (f_*X_1,...,f_*X_n,\xi).
\end{equation}
We have the following Gauss formula
\begin{equation}
\tilde\nabla_X(f_*Y)=f_*(\nabla_XY)+g(X,Y)\xi,
\end{equation}
for vector fields $X,Y$ on $M$. It is known that $\nabla$ is a
torsion-free connection  and $g$ is a symmetric bilinear form on
$M$. The connection $\nabla$ is called the induced connection and
$g$ - the second fundamental form of $f$. The conformal class of $g$
is independent of the choice of a transversal vector field $\xi$. A
hypersurface is called non-degenerate if $g$ is non-degenerate at
each point of $M$.

From now on we shall consider only non-degenerate hypersurfaces.
Hence the second fundamental form $g$ on $M$ is a metric tensor
field (maybe indefinite). The induced volume form on $M$ is, in
general, different than the volume  form $\nu_g$ determined by $g$.
If $\nu_\xi=\nu_g$ we say that the apolarity condition is satisfied.
 A transversal vector field $\xi$ is called equiaffine if
$\nabla\nu_\xi=0$ on $M$. If a transversal vector field is
equiaffine then $\nabla g$ as a cubic form (that is, $\nabla
g(X,Y,Z)=(\nabla _Zg)(Y,Z)$) is totally symmetric.

 The following theorem is central in the
classical affine differential geometry, see \cite{LSZ},\cite{NS}.

\begin{thm}
Let $f:M\to \R^{n+1}$ be a non-degenerate hypersurface. There exists
a unique equiaffine transversal vector field $\xi$ such that
$\nu_\xi=\nu_g$.
\end{thm}

The unique transversal vector field is called the affine normal
vector field. If the  affine lines determined by the affine normal
vector field meet at one point or are parallel then the hypersurface
is called an affine sphere - proper in the first case and improper
in the second one.

Affine spheres are also  described by the affine shape operator.
Namely, let $\xi$ be the affine normal vector field for a
hypersurface immersion $f$. By differentiating $\xi$  relative to
$\tilde \nabla$ we get
\begin{equation}
\tilde\nabla_X\xi = -f_*(SX)
\end{equation}
for some $(1,1)$-tensor field $S$ on $M$. The tensor field $S$ is
called the affine shape operator. The fact that $f$ is an affine
sphere is equivalent to the condition $S=\lambda\,\id$, where
$\lambda $ is a real number, non-zero for a proper sphere and zero
for an improper sphere.

 We have the following
fundamental theorem  for equiaffine hypersurfaces, see \cite{LSZ},
\cite{NS}.
\begin{thm}
Let $\nabla$ be a torsion-free connection on a simply connected
manifold $M$, $g$ be a symmetric bilinear non-degenerate form and
$S$ is a $(1,1)$-tensor field on $M$ such that the following
fundamental equations are satisfied:
\begin{equation}\label{Gauss}
R(X,Y)Z= g(Y,Z)SX-g(X,Z)SY\ \ \ \  \ \ \ -\ Gauss,
\end{equation}
\begin{equation}\label{Ricci}
g(X,SY)=g(SX,Y)\ \ \ \  \ \ \ -\ Ricci,
\end{equation}
\begin{equation}\label{CodazziI}
\nabla g(X,Y,Z)=\nabla g(Y,X,Z)\ \ \ \  \ \ \ -\ I\ Codazzi,
\end{equation}
\begin{equation}
\nabla S(X,Y)=\nabla S(Y,X)\ \ \ \  \ \ \ -\ II \  Codazzi
\end{equation}
for every $X,Y,Z\in T_x M$, $x\in M$, where $R$ is the curvature
tensor for $\nabla$. There is an immersion $f:M\to\R^{n+1}$ and its
equiaffine  transversal vector field $\xi$ such that $\nabla$, $g$,
$S$ are the induced connection, the second fundamental form and the
shape operator for the immersion $f$ equipped with the transversal
vector field $\xi$. The immersion is unique up to an equiaffine
transformation of $\R^{n+1}$. If moreover $\nabla\nu_g=0$ then $\xi$
is  the affine normal (up to a constant) for $f$.
\end{thm}
 For an affine sphere with $S=\lambda\,\id$ the
fundamental equations reduce to the two equations
\begin{equation}\label{fundamental_for_spheres}
\begin{array}{rcl}
&&R(X,Y)Z=\lambda\{g(Y,Z)X-g(X,Z)Y\},\\
&&(\nabla g)(X,Y,Z)=(\nabla g) (Y,X,Z).\\
\end{array}
\end{equation}

As a consequence of the fundamental theorem we have

\begin{corollary}
Let $M$ be a manifold equipped with a metric tensor field $g$, a
torsion-free connection $\nabla$ such that the  equations
{\rm{(\ref{fundamental_for_spheres})}} are satisfied for some
constant real number $\lambda$ and $\nabla \nu_g=0$. For each point
$x$ of $M$ there is a neighborhood $U$ of $x$ and an immersion
$f:U\to \R^{n+1}$ which is an affine sphere whose shape operator is
equal to $\lambda\id$.
\end{corollary}

\proof It is sufficient to define $S=\lambda\, \id$.\koniec

From now on we shall deal with the 2-dimensional case. For a
connection on a 2-dimensional manifold the curvature tensor is
determined by its Ricci tensor. Namely we have
\begin{equation}
R(X,Y)Z=\Ric (Y,Z)X-\Ric(X,Z)Y
\end{equation}
for any vectors $X,Y,Z\in T_xM$, $x\in M$. Therefore the Gauss
equation for a 2-dimensional sphere is equivalent to the condition
$\Ric =\lambda g$.  Hence we have

\begin{corollary}
Let a metric tenor $g$ and a torsion-free connection $\nabla$ be
given on a two-dimensional manifold $M$. They can be locally
realized on an affine sphere if and only if the cubic form $\nabla
g$ is totally symmetric, $\Ric =\lambda g$ for some constant real
number $\lambda$ and $\nabla \nu _g=0$.
\end{corollary}
\medskip

\section{Affine connections, volume forms and the Ricci tensor}
All connections considered in this paper are torsion-free. For any
connection $\nabla$ and a metric tensor field $g$ we denote by $K$
the difference tensor, that is, $K_XY=\nabla_XY-\hat\nabla_XY$,
where $\hat\nabla$ the Levi-Civita connection of $g$. Set
$K(X,Y)=K_XY$. Since both connections $\nabla$ and $\hat\nabla$ are
without torsion, $K$ is symmetric for $X,Y$. The cubic form $\nabla
g$ is symmetric  if and only if $K$ is symmetric relative to $g$,
i.e. $g(K(X,Y),Z)=g(K(X,Z),Y)$ for every $X,Y,Z$. Indeed, we have
$$(\nabla_Xg)(Y,Z)= (K_Xg)(Y,Z)=-g(K(X,Y),Z) -g(Y,K(X,Z)).$$
Assume that the cubic form $\nabla g$ is symmetric. Since
\begin{equation}\label{nabla_Xnu_g}2\nabla_X\nu _g=\tr _g (\nabla _Xg)(\cdot,\cdot
)\nu _g,\end{equation} the condition $ \nabla \gnu =0$
 is equivalent to the condition
\begin{equation}\label{trace-free}
\tr _g(\nabla _Xg)(\cdot,\cdot )=0
\end{equation}
for every $X\in TM$. Since
\begin{equation}
\nabla g(X,Y,Z)=-2g(K(X,Y),Z),
\end{equation}
we have $\nabla\nu_g=0$ if and only if and only if $\tr K_X=0$ for
every $X$.

Recall that the divergence of $X\in\mathcal X(M)$ relative to a
connection $\nabla$ is defined as $\div ^\nabla X=\tr \{Y\to
\nabla_YX\}$.
 Recall also that for a torsion-free connection $\nabla$ there is
(locally) a volume form $\nu$ such that $\nabla\nu =0$ if and only
if the Ricci tensor of $\nabla$ is symmetric. A pair $(\nabla, \nu)$
is then called an equiaffine structure. For a fixed coordinate
system we have
\begin{eqnarray*}
&&(\nabla_{\partial_i}\nu)(\partial_1,...,\partial
_n)=\partial_i(\nu(\partial_1,...,\partial_n))
-\nu(\nabla_{\partial_i}\partial_1,...,\partial_n)-...-\nu(\partial_1,...,\nabla_{\partial_i}\partial_n)\\
&&\ \ \ \ \ \ \ \ \ \ \ \ \ \ \ \ \ \ \ \ \
=\partial_i(\nu(\partial_1,...,\partial_n))
-\nu(\nabla_{\partial_1}\partial_i,...,\partial_n)-...-\nu(\partial_1,...,\nabla_{\partial_n}\partial_i)\\
&& \ \ \ \  \ \ \ \  \ \ \ \ \ \ \ \ \ \ \ \ \ =
\partial_i(\nu(\partial_1,...,\partial_n)) -(\div^\nabla\partial_i)\nu(\partial_1,...,\partial_n).
\end{eqnarray*}
It follows that $\nabla\nu=0$ if and only if $
\partial_i({\rm ln}(\nu(\partial_1,...,\partial_n)))=\div^\nabla
\partial_i$.
Assume now that $\nu _g$ is a volume form determined by a metric
tensor $g$ (not necessarily positive definite). Then
$$ \nu_g(\partial_1,...,\partial _n)^2= G,$$
where $G=| \det[g_{kl}]|$ and  $g_{kl}=g(\partial_k,\partial_l)$ for
$1\le k,l\le n$. Hence $ \nabla\nu _g=0$ if and only if
$\div^\nabla
\partial_i=({\rm ln\,}G)_i/2$ for $i=1,...,n$.

 Let $\Gamma_{ij}$ be the Christoffel symbols of  a connection $\nabla$.
 In general, we have the following formula for the Ricci tensor
 $\Ric $ of $\nabla$
\begin{equation}
\Ric (\pai,\paj)=\sum_{k=1}^n(\Gamma _{ij}^k)_k -(D_j)_i +\Lambda
_{ij}\end{equation} where $D_i=\div ^\nabla \partial _i$ and
\begin{equation}\Lambda
_{ij}=\sum_{k,l=1}^n[\Gamma^l_{kj}\Gamma^k_{il}-\Gamma^l_{ij}\Gamma^k_{kl}].\end{equation}
Since the connection is torsion-free, $\Lambda _{ij}$ is symmetric
for $i$, $j$. From now on we assume that $M$ is $2$-dimensional and
$\Ric$ is symmetric.
 Let $r_{ij}$ be the components of the Ricci tensor in a coordinate
system. Then we have
\begin{equation}\label{Lambda-2-wym}
\begin{array}{rcl}
&&\Lambda
_{11}=\Gamma^2_{11}\Gamma^1_{12}+(\Gamma^2_{12})^2-\Gamma^1_{11}\Gamma^2_{12}-\Gamma^2_{11}\Gamma^2_{22},\\
&&\Lambda _{12}=\Gamma ^1_{22}\Gamma^2_{11}-\Gamma
^1_{12}\Gamma^2_{12},\\
&&\Lambda_{22}=(\Gamma^1_{21})^2 +\Gamma
^2_{12}\Gamma^1_{22}-\Gamma^1_{22}\Gamma^1_{11}-\Gamma^2_{22}\Gamma
^1_{12},
\end{array}
\end{equation}
\begin{equation}\label{CK-2wym}
\begin{array}{rcl}
&&(\Gamma ^2_{12})_1=(\Gamma^2_{11})_2-\Lambda_{11}-r_{11},\\
&&(\Gamma^1_{12})_1=-(\Gamma^2_{12})_2 +\Lambda _{12}+(D_1)_2+r_{12},\\
&&(\Gamma^1_{22})_1=-(\Gamma^2_{22})_2 +\Lambda
_{22}+(D_2)_2+r_{22},
\end{array}
\end{equation}

\begin{equation}\label{D-2-wym}
\begin{array}{rcl}
&&D_1=\Gamma^1_{11}+\Gamma ^2_{21},\\
&&D_2=\Gamma ^1_{12}+\Gamma ^2_{22}.
\end{array}
\end{equation}

Let $g$ be a metric tensor field on  $M$  and $\hat\nabla$ its
Levi-Civita connection. Let $(x^1,x^2)$ be an isothermal coordinate
system for $g$, that is, $g_{22}=\varepsilon g_{11}$,
$g_{11}=e^\varphi$ and $g_{12}=0$ for some function $\varphi$ and
$\varepsilon =\pm$ depending on the signature of $g$. For an
isothermal coordinate system we have
\begin{equation}
\begin{array}{rcl}
&&\hat\nabla _{\partial_1}{\partial_1}=\frac{\varphi_1}{2}\partial
_1-\varepsilon \frac{\varphi_2}{2}\partial _2,\\
&&\hat\nabla _{\partial_1}{\partial_2}=\frac{\varphi_2}{2}\partial
_1+\frac{\varphi_1}{2}\partial _2,\\
&&\hat\nabla
_{\partial_2}{\partial_2}=-\varepsilon\frac{\varphi_1}{2}\partial
_1+\frac{\varphi_2}{2}\partial _2.
\end{array}
\end{equation}
The curvature $\kappa$ of $g$ is given by
\begin{equation}\label{curvature}
\kappa =-\frac{\Delta\varphi}{2}.
\end{equation}
Recall that for any function $\alpha$ we have
\begin{equation}\label{laplacian}
\Delta \alpha=\frac{\alpha_{11}+\varepsilon \alpha_{22}}{e^\varphi}.
\end{equation}

Assume  $\nabla$ is a  connection on $M$ such that the cubic form
$\nabla g$ is totally symmetric. The difference tensor
$K=\nabla-\hat\nabla$ is  symmetric relative to $g$. For the
isothermal coordinate system the symmetry is equivalent to the
conditions
\begin{equation}\label{K-symmetry}
K^2_{21}=\varepsilon K^1_{22}, \ \ \  K^1_{12}=\varepsilon
K^2_{11}.\end{equation}
Assume moreover that  $\nabla\nu_g=0$. Then
$\div ^\nabla\partial_1=\varphi_1$ and $\div^\nabla
\partial_2=\varphi _2$. Since
\begin{eqnarray*}
&&\varphi_1=\div ^\nabla\partial_1=\hat\Gamma^1_{11}+K^1_{11}+\hat\Gamma^2_{21}+K^2_{21},\\
&&\varphi_2=\div
^\nabla\partial_2=\hat\Gamma^1_{12}+K^1_{12}+\hat\Gamma^2_{22}+K^2_{22},
\end{eqnarray*}
we get
\begin{equation}\label{K-apolarity}
K^1_{11} =-K^1_{22},\ \ \  \ K^2_{11}=-K^2_{22}.
\end{equation}
We see that among the functions $K^i_{jk}$, $i,j,k=1,2$, only two
functions are independent. We choose the functions $K^1_{12}$ and
$K^2_{21}$. Using (\ref{Lambda-2-wym}), (\ref{K-apolarity}) and
(\ref{K-symmetry}), by a straightforward computation, one gets
\begin{equation}
\begin{array}{rcl}\label{Lambda}
&&\Lambda_{11}=\varphi _1K^2_{21} -\varepsilon
\varphi_2K^1_{12}+2[(K^2_{21})^2 +\varepsilon(K^1_{12})^2],\\
&&\Lambda_{12}=-\varphi_1K^1_{12}-\varphi_2K^2_{21},\\
&&\Lambda _{22}=\varphi_2 K^1_{12} -\varepsilon\varphi_1K^2_{21}
+2[(K^1_{12})^2+\varepsilon(K^2_{21})^2] .
\end{array}
\end{equation}
We have
\begin{equation}
g(K,K)= \frac{4}{e^\varphi}[\varepsilon(K^1_{12})^2+(K^2_{21})^2].
\end{equation}
In the theory of affine surfaces, the function $\frac{1}{2}g(K,K)$
is called the Pick invariant and it is usually denoted by $J$.

Formulas (\ref{CK-2wym}) for the Ricci tensor $\Ric=r_{ij}$ of the
connection $\nabla$, after using (\ref{Lambda}), receive the
following form:
\begin{equation}\label{Ricci_K}
\begin{array}{rcl}
&&(K^2_{21})_1=-\frac{\varphi_{11}+\varepsilon\varphi_{22}}{2}
+\varepsilon
(K^1_{12})_2\\
&&\ \ \ \ \ \ \ \ \ \ \ \ \ \ \ -\varphi_1
K^2_{21}+\varepsilon\varphi_2K^1_{12}
-2[\varepsilon(K^1_{12})^2+(K^2_{21}) ^2]-r_{11},\\
&&(K^1_{12})_1=-(K^2_{21})_2-\varphi_1K^1_{12}-\varphi_2K^2_{21}+r_{12},\\
&&(K^2_{21})_1=
\frac{\varphi_{11}+\varepsilon\varphi_{22}}{2}+ \varepsilon(K^1_{12})_2\\
&&\ \ \ \ \ \ \  \ \ \ \ \  \ \ \
+\varepsilon\varphi_2K^1_{12}-\varphi_1K^2_{21}
+2\varepsilon[\varepsilon (K^1_{12})^2+(K^2_{21})^2]+\varepsilon
r_{22}.
\end{array}
\end{equation}
When adding and subtracting the first and the last equalities we get
the following system of equalities equivalent to (\ref{Ricci_K})
\begin{equation}\label{uklad}
\begin{array}{rcl}
&&r_{11}+\varepsilon
r_{22}=-(\varphi_{11}+\varepsilon\varphi_{22})-4[\varepsilon(K^1_{12})^2+(K^2_{21})^2],\\
&&(K^2_{21})_1=\varepsilon (K^1_{12})_2+
\varepsilon\varphi_2K^1_{12}-\varphi_1K^2_{21}+\frac{-r_{11}+\varepsilon
r_{22}}{2},\\
&&(K^1_{12})_1=-(K^2_{21})_2-\varphi_1K^1_{12}-\varphi_2K^2_{21}+r_{12}.
\end{array}
\end{equation}
The first equality is the affine theorema egregium, that is, the
equality
\begin{equation}\label{theorema_egregium}\rho=\kappa-J,\end{equation}
 where $\rho=(\tr _g\Ric)/2$.
\medskip

\section{Proof of Theorem \ref{main}}
 Direct proofs of necessity of the condition  (\ref{warunek}) can be found, for instance, in \cite{NS} or \cite{LSZ}.
The proof below  also gives the necessity, but we focus on the
sufficiency of the condition. Assume that for a given $g$ and some
constant $\lambda$ the equality (\ref{warunek}) is satisfied. By the
fundamental theorem we are looking for a connection $\nabla$ such
that the cubic form $\nabla g$ is symmetric, the Ricci tensor of
$\nabla$ is equal to $\lambda g$ and $\nabla \nu_g=0$. Let
$\nabla=\hat\nabla+K$.  As in the previous section we will carry
considerations for a fixed   isothermal coordinate system for $g$.
Instead of looking for a connection $\nabla$ we will look for a
tensor field $K$ satisfying appropriate symmetry conditions and the
system of differential equations (\ref{uklad}), where $\lambda=
\frac{r_{11}+\varepsilon r_{22}}{2e^\varphi}$. It will turn out that
(\ref{warunek}) is the integrability condition for the system.

Since $g(K,K)$ should be non-zero, the tensor $K$ should be
non-zero. Suppose that $K^2_{21}\ne 0$. Set
\begin{equation}\label{definicje_L_f}
L=\varepsilon(K^1_{12})^2 +(K^2_{21})^2, \ \  l=
\frac{K^1_{12}}{K^2_{21}}.
\end{equation}

The two functions $L$ and $l$ determine the difference tensor $K$ up
to sign. This is sufficient for our consideration because affine
spheres go in pairs. More precisely, if  $g$, $\nabla=\hat\nabla
+K$, $S=\lambda \,\id$ constitute the induced structure on an affine
sphere then $g$, $\overline\nabla=\hat\nabla -K$, $S=\lambda \,\id$
form the induced structure on another affine sphere. It follows from
the fact that for an affine sphere $R=\overline R$, where $R$ and
$\overline R$ are the curvature tensors for $\nabla$ and $\overline
\nabla$.

We now have
\begin{equation}\label{K221}
(K^2_{21}) ^2=\frac{L}{\varepsilon l^2+1}.\end{equation}
 Set
$r=\lambda g$.
 The system (\ref{uklad}) now becomes
\begin{equation}\label{uklad1}
\begin{array}{rcl}
&&2\frac{L}{e^\varphi}=\kappa -\lambda,\\
&&(K^2_{21})_1=\varepsilon (K^1_{12})_2+
\varepsilon\varphi_2K^1_{12}-\varphi_1K^2_{21},\\
&&(K^1_{12})_1=-(K^2_{21})_2-\varphi_1K^1_{12}-\varphi_2K^2_{21}.
\end{array}
\end{equation}
Note that the denominator in (\ref{K221}) is different than $0$ if
and only if $L\ne 0$, that is, by the first equality from
(\ref{uklad1}), if and only if $\kappa-\lambda=2L/e^\varphi\ne 0$.
Of course, we should define $L=
\frac{(\kappa-\lambda)e^\varphi}{2}$. The function $L$ is now given.

We want to prove that the system of the last two equations from
(\ref{uklad1}) relative to unknown functions $K^1_{12}$, $K^2_{21}$
has a solution. Using (\ref{definicje_L_f}) and(\ref{uklad1}) we
obtain
\begin{equation}\label{L_1}
\begin{array}{rcl}
&& L_1=2\{\varepsilon(K^1_{12})_1K^1_{12}+(K^2_{21})_1K^2_{21}\}\\
&&\ \ \ \ \ \ \ \ \ \  =2\{
\varepsilon[-(K^2_{21})_2-\varphi_1K^1_{12}-\varphi_2K^2_{21}]K^1_{12}\\
&&\ \ \ \ \ \ \ \ \ \ \ \ \ \ \ \ \  +[\varepsilon (K^1_{12})_2+
\varepsilon\varphi_2K^1_{12}-\varphi_1K^2_{21}]K^2_{21}\}\\
&&\ \ \ \ \ \ \ \ \ \ \ \ \ \ \ \ \ \ \ \ \ \ \ \ \ =2\{\varepsilon
l_2(K^2_{21})^2-\varphi _1 L\}
\end{array}
\end{equation}
 and consequently
\begin{equation}
\frac{L_1}{2L}=\frac{l_2}{ l^2+\varepsilon}-\varphi_1.
\end{equation}
Similarly, using (\ref{definicje_L_f}) and (\ref{uklad1}), we get
\begin{equation}\label{f_1}
\begin{array}{rcl}
&&l_1=
\frac{(K^1_{12})_1K^2_{21}-(K^2_{21})_1K^1_{12}}{(K^2_{21})^2}\\
&&\ \ \
 \ \ \ \ =\frac{[-(K^2_{21})_2-\varphi_1K^1_{12}-\varphi_2K^2_{21}]K^2_{21}-[\varepsilon
(K^1_{12})_2+
[\varepsilon\varphi_2K^1_{12}-\varphi_1K^2_{21}]K^1_{12}]}{(K^2_{21})^2}\\
&&\ \ \ \ \  \ \ \ \ =
\frac{-\frac{L_2}{2}-{\varphi_2L}}{(K^2_{21})^2}
\end{array}
\end{equation}
and consequently (by (\ref{K221}))
\begin{equation}
\frac{l_1}{l^2+\varepsilon}=-\varepsilon\frac{L_2}{2L}-\varepsilon\varphi_2.
\end{equation}
We have got the following system of differential equations relative
to $l$
\begin{equation}\label{main_equation}
\begin{array}{rcl}
&&\frac{l_1}{l^2+\varepsilon}=-\varepsilon(\frac{1}{2}\ln |L|+\varphi)_2,\\
&&\frac{l_2}{l^2+\varepsilon}=(\frac{1}{2}\ln |L|+\varphi)_1.
\end{array}
\end{equation}
We have $\frac{l_1}{l^2+\varepsilon}=F_1$ and $\frac{l_2}{
l^2+\varepsilon}=F_2$, where $F=arc\, tg\, l$ or $F= \ln
|\frac{l-1}{l+1}|$ depending on whether $\varepsilon=1$ or $-1$. The
integrability  condition for (\ref{main_equation}) is
$$-\varepsilon(\frac{1}{2}\ln |L|+\varphi)_{22}=(\frac{1}{2}\ln
|L|+\varphi)_{11},$$ which is equivalent to
$$\frac{(\ln |L|)_{11}+\varepsilon(\ln |L|)_{22}}{4e^\varphi}=
-\frac{\varphi_{11}+\varepsilon\varphi_{22}}{2e^\varphi},
$$
that is,
$$\Delta \ln|L|=4\kappa.$$
The last condition is equivalent to

\begin{equation}
\Delta\ln\frac {2|L|}{e^\varphi}=6\kappa,
\end{equation}
which is the desired condition (\ref{warunek}).

Thus the system (\ref{main_equation}) has a solution $l$ around any
fixed point $p\in M$. As the initial condition we  take
$l(p)=\beta$, where $\beta$ can be any real number if $\varepsilon
=1$ and $\beta^2\ne 1$ if $\varepsilon =-1$ (which  again is the
condition $J\ne 0$). From (\ref{K221}) and (\ref{definicje_L_f}) we
get $K^1_{12}$ and $K^2_{21}$. The other components of the tensor
$K$ are defined by using formulas (\ref{K-symmetry}) and
(\ref{K-apolarity}).

 We shall now  check that the obtained functions $K^1_{12}$,
$K^2_{21}$ satisfy the last two equations of (\ref{uklad1}) if
$L=\varepsilon (K^1_{12})^2+(K^2_{21})^2$. Using also the definition
of $l$ one gets
\begin{equation}\label{algebraic_equality}
\begin{array}{rcl}
&&L_1=2\{\varepsilon(K^1_{12})_1K^1_{12}+(K^2_{21})_1K^2_{21}\}\\
&&(K^2_{21})^2l_1=(K^1_{12})_1K^2_{21}-(K^2_{21})_1K^1_{12}.
\end{array}
\end{equation}
We have already checked that  if we substitute the quantities
$(K^1_{12})_1$, $(K^2_{21})_1$ by the right hand sides of the last
two formulas of (\ref{uklad1}), then we get the equalities. We shall
now regard (\ref{algebraic_equality})  as a system of algebraic
linear equations with unknowns $(K^1_{12})_1$, $(K^2_{21})_1$. The
main determinant of (\ref{algebraic_equality}) is equal to $-2L\ne
0$, hence the system of equations has only one solution. It means
that the functions $K^1_{12}$, $K^2_{21}$ satisfy the last two
equations of (\ref{uklad1}).  It follows (by (\ref{uklad})) that the
Ricci tensor of $\nabla=\hat\nabla+K$ satisfies the conditions
$\Ric(\partial _1,\partial _1)=\varepsilon\Ric(\partial _2,\partial
_2)$ and $\Ric(\partial _1,\partial _2)=0$.  Since the scalar
curvature, say $\rho$, of $\nabla$ satisfies the equality
$J+\rho=\kappa$ and we have $J+\lambda=\kappa$, we see that
$\Ric=\lambda g$. \koniec

\begin{remark}
{\rm It is clear from the above proof that an affine sphere whose
Blaschke metric is prescribed (and satisfies the condition
(\ref{warunek})) is not unique. First of all the function $l$ can be
prescribed at a point $p$.  Different functions $l$ give different
(non-equivalent modulo the affine special  group acting on $\R^3$)
affine spheres. Of course, the difference tensors $K$ and $-K$ also
give two different spheres if $K\ne 0$.}
\end{remark}
\section{Proof of Corollaries \ref{corollary_1} and
\ref{corollary_2}}

 Assume first that  $g$ satisfies  (\ref{wstep1}) and
 $(x^1,x^2)$ is an isothermal coordinate system for $g$ as in the previous sections. Denote by
 $\Delta_0$ the Laplacian for the coordinate system, that is,
 $\Delta_0\alpha=\alpha_{11}+\varepsilon \alpha_{22}$ for a function $\alpha$. This is
 the Laplacian for the flat metric tensor field $g_0$ adapted to the
 coordinate system, that is, $g_0(\partial_1,\partial_1)=1$,
 $g_0(\partial_1,\partial_2)=0$,  $g_0(\partial_2,\partial_2)=\varepsilon$,
 where $\varepsilon=\pm 1$.
The equality (\ref{wstep1}) is equivalent to the equality
\begin{equation}
\Delta_0[\ln|\kappa -\lambda|+3\varphi]=0.\end{equation} Hence
\begin{equation}\label{corollary_4}
\ln|\kappa -\lambda|+3\varphi=h\end{equation}
 for some $\Delta_0$-harmonic function  $h$. Thus
 $\kappa-\lambda =\pm e^{h-3\varphi}$. Since $\kappa
 =-\frac{\Delta_0\varphi}{2e^\varphi}$ (by (\ref{curvature}) and (\ref{laplacian})), the last equality is
 equivalent to the equality
 \begin{equation}\label{corollary_3}
-\frac{\Delta_0 c}{2e^c}=\mp 2-\lambda e^{\frac{h-3c}{2}},
 \end{equation}
 where $c=h-2\varphi$.
If $\lambda=0$ then the last equality becomes
\begin{equation}
-\frac{\Delta_0 c }{2e^c}=\pm 2.
\end{equation}
It means, by (\ref{curvature}) and (\ref{laplacian}), that $c$ is
such a function that that the metric $\overline g=e^cg_0$ has
constant Gaussian curvature $\pm 2$. Note that the metric $\overline
g$ is defined only locally, that is, an a domain of an isothermal
coordinate map.
\medskip

{\it{Proof of Corollary {\rm \ref{corollary_1}}}}. Let $f:M\to \R^3$
be an improper locally strongly convex affine sphere with
nowhere-vanishing Pick invariant, equivalently, with nowhere
vanishing curvature $\kappa$ of the Blaschke metric $g$. Hence the
above consideration can be applied. The atlas of isothermal
coordinates is analytic. The Levi-Civita connection for $\overline
g$ is locally symmetric and therefore analytic. It follows that its
curvature tensor is analytic and so is its Ricci tensor. The Ricci
tensor is equal to $\pm 2\overline g$. Hence $\overline g$ is
analytic and $c$ is analytic. If $g$ is definite then a $\overline
g$-harmonic function is analytic. It follows that $\varphi
=\frac{h-c}{2}$ is analytic and consequently
$g=e^{\frac{h-3c}{2}}\overline g=e^{\frac{h-c}{2}}g_0$ is analytic.

We have proved that  the Blaschke metric $g$ is  analytic. An
improper affine sphere can be locally regarded as a graph of some
function and its affine normal is a constant vector. More precisely,
$f$ is, up to affine transformations of $\R^3$, given locally by
$$\R^2 \supset U\ni (x^1,x^2)\to (x^1,x^2, \Psi (x^1,x^2))\in R^3,$$
where $\Psi$ is a smooth function and its affine normal is equal to
$(0,0,1)$.  We have
$$g(\partial_i,\partial _j)=\Psi_{ij}.$$
Since $g$ is analytic, so is $\Psi$. The proof of Corollary
\ref{corollary_1}   is completed. \koniec

\medskip
{\it Proof of  Corollary {\rm \ref{corollary_2}}}. Assume now that
$h$ is an arbitrary harmonic function on some open set $U\subset
\R^2$ and $c$ is such a function on $U$ that $g=e^cg_0$ has constant
Gaussian curvature $2$ or $-2$. Of course, such a function exists
because there exist metrics of any constant Gaussian curvature. Set
$\varphi =\frac{h-c}{2}$. By the consideration from the beginning of
this section one sees that  the equality (\ref{corollary_4}) is
satisfied for $\kappa =-\frac{\Delta_0\varphi}{2e^\varphi}$.\koniec

\end{document}